\documentclass[11pt, twoside]{article}
\usepackage{amsfonts,amssymb,amsmath,amsthm}
\usepackage{authblk}
\usepackage{graphicx}
\usepackage[all]{xy}
\usepackage{tikz}
\usetikzlibrary {arrows.meta}
\usepackage{changepage}
\usepackage{psfrag,xmpmulti,amscd,color,pstricks, import,enumerate}
\usepackage[normalem]{ulem}
 \divide\oddsidemargin by 2
\newtheorem{thm}{Theorem}[section]

\newtheorem{lem}[thm]{Lemma}

\theoremstyle{definition}
\newtheorem{defn}[thm]{Definition}

\newtheorem*{rem*}{Remark}
\newtheorem*{rems*}{Remarks}

\newtheorem*{ex*}{Example}

\numberwithin{equation}{section}

\definecolor{OrangeRed}{cmyk}{0,0.6,1,0}            % half magenta only, full yellow
\definecolor{DarkBlue}{cmyk}{1,1,0,0.20}
\definecolor{DarkGreen}{cmyk}{1,0,0.6,0.2}
\definecolor{myblue}{rgb}{0.66,0.78,1.00}
\definecolor{Violet}{cmyk}{0.79,0.88,0,0}
\definecolor{Lavender}{cmyk}{0,0.48,0,0}
\definecolor{purpleheart}{rgb}{0.41, 0.21, 0.61}
\definecolor{brick}{cmyk}{0,0.8,0.3,0.5}

\newcommand{\C}{{\mathbb C}}

\newcommand{\D}{{\mathbb D}}
\newcommand{\Hyp}{{\mathbb H}}

\newcommand{\N}{{\mathbb N}}

\newcommand{\R}{{\mathbb R}}

\newcommand{\Z}{{\mathbb Z}}

\newcommand{\ov}{\overline}

\renewcommand{\epsilon}{\varepsilon}
\renewcommand{\phi}{\varphi}

\newcommand{\tef}{transcendental entire function}

\vspace{5cm}

\title{Quasiconformal folding - a review of `Models for the Eremenko--Lyubich class'}

\author[1]{Philip J. Rippon}

\affil{\small School of Mathematics and Statistics, The Open University, Milton Keynes MK7 6AA\\phil.rippon@open.ac.uk}

\date{\today}

\begin{document}

\maketitle
\begin{abstract}
We review the paper `Models for the Eremenko--Lyubich class' by Chris Bishop, which appeared in the Journal of the London Mathematical Society in 2015, the first of two LMS papers for which Bishop received the LMS Senior Berwick Prize in 2024. This is one of a series of papers where Bishop introduced his new technique of quasiconformal folding, which transformed the range of approaches to constructing examples of transcendental entire functions leading to the resolution of many open problems, in particular in transcendental dynamics.

\end{abstract}

\section{Background and main result}

The papers \cite{BishopJLMS, BishopActa, BishopPLMS} by Chris Bishop appeared within just a few years, each containing versions of a new technique he had introduced, called quasiconformal folding. This technique can be used to construct transcendental entire functions, and other holomorphic functions, that lie in various different important classes and satisfy a large variety of different conditions. This review will focus on the results and techniques in \cite{BishopJLMS}, which Bishop describes as a `gentle introduction' to quasiconformal folding, and will briefly discuss Bishop's related papers, applications of his results, and further developments.

The origin of the paper \cite{BishopJLMS} can be traced back to a conversation that took place at Stony Brook in March 2011, when Alex Eremenko was visiting Misha Lyubich, with whom he had collaborated extensively in earlier years on problems in transcendental dynamics; see~\cite{EandL84, EandL87,EandL}. Eremenko had a question for Chris Bishop during a conversation that started in a stairwell. The question amounted to asking whether every connected compact set in the complex plane can be approximated arbitrarily closely by the critical points of a polynomial with only two critical values.

Bishop's track record of solving technically extremely hard problems in complex analysis, involving (amongst many other things) harmonic measure and quasiconformal mapping, meant that he was ideally prepared for such a problem. Moreover, his solution of Eremenko's problem in \cite{BishopInv} led him to the ingenious and complicated technique that he called quasiconformal folding. This technique, along with related ideas, enabled Bishop to solve in a remarkably short period of time several major open questions about entire functions, in particular ones in transcendental dynamics, a field where he had not previously worked! It is no exaggeration to say that his work on quasiconformal folding represents the most transformative introduction of a mathematical technique in transcendental dynamics this century.

Before describing Bishop's work in more detail, we mention that a good introduction to transcendental dynamics can be found in the survey article~\cite{bergweiler93} and an encyclopedic account of the \textit{escaping set}, $I(f)=\{z\in\C:f^n(z)\to\infty\;\text{as }n\to\infty\}$, where~$f$ is a transcendental entire or meromorphic function, is in the recent survey article~\cite{BergRemp25}.

Two key classes of \tef s are the Speiser class and the Eremenko--Lyubich class, which both consist of functions that have a restricted range of inverse function singularities. The dynamical behaviour of functions in these classes turns out to be more tractable than that of general entire functions, and hence they have been the subject of intense research.

For a {\tef}~$f$ a \textit{singularity} of $f^{-1}$ is any finite point in some neighbourhood of which some inverse branch of~$f$ cannot be analytically continued. For example, 0 is a singularity of $f^{-1}$ for $f(z)=e^z$; indeed, no inverse branch of~$f$ can be analytically continued in any neighbourhood of 0. It can be shown that for any transcendental entire function~$f$ the singularities of $f^{-1}$ are precisely the \textit{critical values} of~$f$, defined as the images of critical points of $f$, and the \textit{asymptotic values} of~$f$, defined as the finite limit values of~$f$ along paths in $\C$ tending to~$\infty$. The set of singularities of~$f^{-1}$ is denoted by ${\rm sing}( f^{-1})$, and its closure $S(f)$ is called the \textit{singular set}. Thus the complement of $S(f)$ is exactly the set of points that have a neighbourhood~$U$ in which $f: f^{-1}(U) \to U$ is unbranched.

The \textit{Speiser class} is
\[
\mathcal S:=\{f: f \text{ is a {\tef} and } {\rm sing}(f^{-1}) \text{ is finite}\},
\]
and the \textit{Eremenko--Lyubich class} is
\[
\mathcal B:=\{f: f \text{ is a {\tef} and } {\rm sing}(f^{-1}) \text{ is bounded}\}.
\]
Clearly $\mathcal S \subset \mathcal B$. Simple examples of functions in the class $\mathcal S$ are $f(z)=e^z$ and $f(z)=\sin z$, whereas the entire function $f(z)=\sin z/z$ lies in the class~$\mathcal B$ but not in the class~$\mathcal S$ (since it has infinitely many critical values).

Functions with just a finite set of inverse function singularities were studied almost one hundred years ago by Speiser as part of a programme to understand Riemann surfaces, and they came to the fore in transcendental dynamics when Eremenko and Lyubich \cite{EandL84, EandL}, and also Goldberg and Keen \cite{GK}, showed that functions in the class $\mathcal S$ do not have \textit{wandering domains}, components of the Fatou set (the stable set under iteration of the function) whose orbits do not end in a cycle. Their proofs adapted techniques from the theory of quasiconformal mappings used by Sullivan to prove his celebrated no wandering domains theorem for rational maps. It had been known for some time that {\tef}s can have wandering domains, either simply connected or multiply connected (see \cite{Baker76} and \cite{bergweiler93} for many examples) and their detailed study is one of the main themes in contemporary transcendental dynamics; see for example, \cite{lasse-helena,  brs, bfjk, BT, BEFRS, EGP, MRW}.

At about the same time, Eremenko and Lyubich \cite{EandL} introduced the class $\mathcal B$ (later named the Eremenko--Lyubich class by Rempe) and showed that functions in this larger class never have Fatou components in the escaping set $I(f)$. But this left open the intriguing question whether some functions in the class~$\mathcal B$ could have an oscillating orbit of wandering domains with some subsequence accumulating in a bounded part of the plane.

In another paper \cite{EandL87}, which is reviewed in this anniversary volume~\cite{FagePard}, Eremenko and Lyubich had used approximation theory, for the first time in complex dynamics, to construct the first example of a {\tef} with such an oscillating orbit of wandering domains, which was not however in the class $\mathcal B$. This intriguing question about the possible existence of wandering domains in the class~$\mathcal B$ was eventually answered by Bishop, as we shall see.

To prove their result about functions in the class~$\mathcal B$, Eremenko and Lyubich established a fundamental property of functions in this class: they showed that if the set of singular values lies in a disc $\{z\in\C : |z| < R\}$, then the preimage of the exterior of the disc, that is, $\{z\in\C : |f(z)| > R\}$, is a disjoint union of unbounded Jordan domains, $\Omega_j$ say, and also that~$f$ acts as a covering map from each component $\Omega_j$ to $\{z\in\C : |z| > R\}$. This covering property implies that the function~$f$ has significant expansion properties in these components, similar to the expansion properties of the exponential function in the right half-plane.

The main result in Bishop's paper \cite{BishopJLMS} reverses this fundamental process, that is, it starts with a set~$\Omega$ which is a disjoint union of unbounded simply connected domains and their associated covering maps, which is to be thought of as a `model', and constructs an entire function in the class~$\mathcal B$ that approximates these covering maps. This process reduces the problem of constructing a function in the class~$\mathcal B$ with desired properties to the usually easier problem of constructing an appropriate model with those properties.

We now give the definition of a model in \cite{BishopJLMS}, which involves conformal (Riemann) mappings to the right half-plane, denoted by $\mathbb{H}_r = \{ x+iy  : x > 0\}$; see Figure~1, borrowed from~\cite{BishopJLMS}.

\begin{defn} Let $I$ be an at most countably infinite index set and let $\Omega = \bigcup_{j \in I} \Omega_j$ be a disjoint union of unbounded simply connected domains $\Omega_j$. For each $j\in I$, let $\tau_j : \Omega_j \to \mathbb{H}_r$ be conformal, and $\tau: \Omega\mapsto \mathbb{H}_r$ equal $\tau_j$ on $\Omega_j$, for $j \in I$. Also, suppose that the following conditions hold:
\begin{enumerate}[(a)]
\item Sequences of components of $\Omega$ accumulate only at $\infty$.\label{condaccum}
\item The boundary of each $\Omega_j$,  $j \in I$, is connected (in $\C$).\label{condboundary}
\item For any sequence $(z_n)_{n\in\N}$ in $\Omega$, if $\tau(z_n)\rightarrow \infty$ as $n\rightarrow\infty$, then $z_n\rightarrow\infty$ as $n\rightarrow\infty$.\label{condtoinf}
\end{enumerate}
Finally, set $F(z) = \exp(\tau(z))$, $z\in\Omega$, which is a covering map from each $\Omega_j$ onto $\C\setminus\overline{\D}$. The pair $(\Omega, F)$ is called a \emph{model} and the domains $\Omega_j$ are called \textit{tracts}.
\end{defn}

\begin{figure}
\centering
\includegraphics[width=0.8\textwidth]{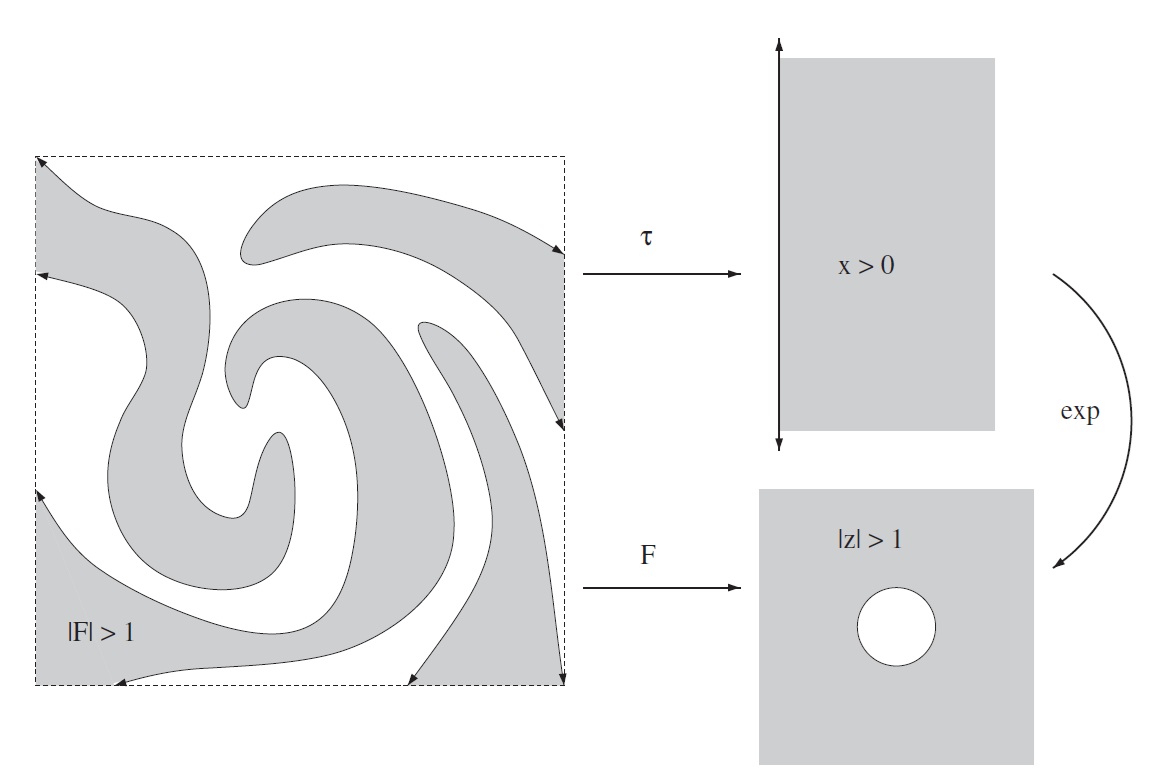}
\caption{\label{model} \small A model consists of tracts with covering maps defined via conformal mappings to the right half-plane}
\end{figure}

Given a model $(\Omega, F)$ and $\rho>0$, we define these subsets of $\Omega$:
\[
\Omega(\rho)=\{z\in\C:|F(z)|>e^\rho\}=\tau^{-1}(\{x+iy:x>\rho\}),\quad \Omega_j(\rho)=\Omega(\rho)\cap \Omega_j.
\]
Here is a version of Bishop's main theorem in \cite{BishopJLMS}; some of the original detailed conclusions of the theorem have been omitted for simplicity of exposition.
\begin{thm}\label{Chris-jlms-main}
Suppose that $(\Omega,F)$ is a model and $0<\rho\le 1$. Then there is an entire function $f\in \mathcal B$ and a quasiconformal mapping $\phi:\C\to \C$ such that $F=f\circ \phi$ on $\Omega (2\rho)$ and
\[
|f\circ\phi(z)|\le e^\rho,\;\;\text{for } z \notin \Omega(\rho),\quad\text{and}\quad |f\circ\phi(z)|\le e^{2\rho},\;\;\text{for } z \notin \Omega(2\rho).
\]
In particular, $S(f)\subset \{z\in \C:|z|\le e^\rho\}$ and the components of $\{z\in\C:|f(z)|>e^\rho\}$ are in 1-to-1 correspondence with the components of $\Omega$ via the mapping~$\phi$.
\end{thm}
To make this result clear, we need to explain the term `quasiconformal mapping'. Briefly, a {\it quasiconformal mapping} $\phi$ is a homeomorphism between plane domains with the property that at almost all points~$z$ of the domain it is differentiable (in the sense that its derivative is a linear self map of $\R^2$) and, intuitively,~$\phi$ maps infinitesimal circles centred at~$z$ to infinitesimal ellipses centred at $\phi(z)$. If the dilatation of these ellipses, that is, the ratio of the length of the major axis of the ellipse to the length of the minor axis, is almost always at most~$K$, then the mapping is called $K$-quasiconformal; see \cite{AhlforsQCmaps} or \cite{BranFage} for more details.

For example, any $\R$-linear mapping on $\C$ is quasiconformal and maps circles to ellipses for which the dilatation is constant. In general, quasiconformal mappings need not be $C^{(1)}$ but they satisfy a global property that is slightly weaker than Lipschitz (see \cite[Chapter~II, §4.2]{LehtVirt}), so Theorem~\ref{Chris-jlms-main} guarantees that prescribed geometric properties of the tracts $\Omega_j$ of the model are reproduced in the corresponding components $\{z\in\C:|f(z)|>e^\rho\}$ of~$f$, also known as tracts.

In Section~2, we outline the main ideas of the proof of Theorem~\ref{Chris-jlms-main} and in Section~3 we describe Bishop's further results on the construction of functions in both classes $\mathcal B$ and $\mathcal S$ using quasiconformal folding, together with some subsequent related developments.

\begin{rem*}
The interested reader is recommended to consult material on Bishop's website at Stony Brook, where he lists comments and corrections for his papers on quasiconformal folding and also includes his previous postgraduate lectures on many topics. In particular, there is a rich and lengthy set of lecture notes on transcendental dynamics, with an introduction describing the history of his invention of the technique of quasiconformal folding \cite{BishNotes}.
\end{rem*}

\section{Outline proof of Theorem~\ref{Chris-jlms-main}}

In this section, we outline Bishop's proof of Theorem~\ref{Chris-jlms-main}. For simplicity of notation, we assume that the open set~$\Omega$ in the model has a single component, so there is no need for an index set~$I$. The proof assumes that the reader has considerable knowledge of Blaschke products, conformal mappings and their distortion properties, quasiconformal mappings, quasiregular mappings, the hyperbolic metric and harmonic measure.

We start with a brief overall summary. The first step is to show  that it is sufficient to consider the case $\rho=1$. Having done that, Bishop introduces a new holomorphic map~$G$, which is the composition of a Riemann mapping from the unbounded simply connected domain $W=\C\setminus \ov{\Omega(1)}$ to the open unit disc~$\D$, followed by a carefully chosen self map of~$\D$ which is an infinite Blaschke product~$B$ whose zeros depend on the geometry of the boundary of $\Omega(1)$. Then, an interpolation takes place across the set $\Omega(1)\setminus \ov{\Omega(2)}$ between the new map~$G$ on~$W$ and the original map~$F$ restricted to $\Omega(2)$, using an ingenious composition of quasiconformal mappings (this is where the folding happens).

This interpolation process results in a so-called quasiregular mapping~$g$ defined on the whole of~$\C$ with uniformly bounded dilatation. A \textit{quasiregular mapping} can be defined as a mapping that is locally quasiconformal except at a discrete set of `critical' points where it behaves locally like a power map $z\mapsto z^d$. The proof is completed by applying the Measurable Riemann Mapping Theorem, arguably the major result in quasiconformal theory, due in various forms to Morrey, Bojarski, and Ahlfors and Bers. This gives an entire function~$f$ and a quasiconformal mapping~$\phi:\C\to \C$ such that
\[
g=f\circ \phi\;\;\text{on }\C \quad \text{and hence}\quad F=f\circ \phi\;\;\text{on }\Omega (2),
\]
as required; for proofs of the Measurable Riemann Mapping Theorem and these properties of quasiregular mappings, see \cite[Sections~1.4 and~1.6]{BranFage}.

Now we give some details of the interpolation process, shown in Figure~2.  Most of the construction takes place in the strip~$S$ in $\Hyp_r$ bounded by the lines $L_1=\{x+iy:x=1\}$ and $L_2=\{x+iy:x=2\}$.
%let $\gamma=\partial W$, and let $\mathcal{J}_\gamma$ be the partition of $\gamma$ that is the image of the partition $\mathcal J$ under $\tau^{-1}$.
Recall that $\Psi$ is a Riemann mapping from $W=\C\setminus \ov{\Omega(1)}$ to the open unit disc~$\D$. Then, by an application of the Schwarz reflection principle, the map~$\Psi$ extends to $\C\setminus \ov{\Omega(2)}$ as a conformal mapping and in particular as a smooth map of $\partial W$ onto $\partial \D$ apart from one point of $\partial\D$, shown as $1$ in Figure~2, which can be thought of as the image of $\infty$.

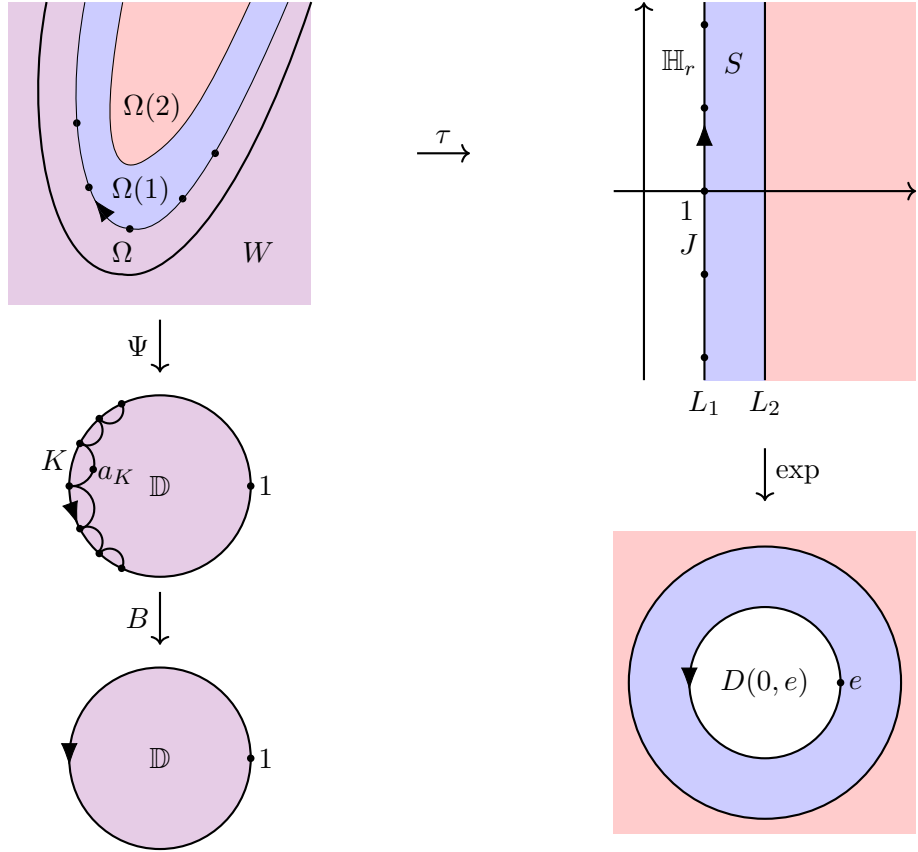
\begin{figure}[hbt!]
\centering
\begin{tikzpicture}[scale=1]
\tikzset{centerdot/.style={circle, fill=black, inner sep=1pt}}

%---------------- LEFT COLUMN: Omega -> D -> D ----------------%

% Omega domain
\begin{scope}[shift={(0,0)}]
  % Outer violet region
  \fill[violet!20] (-2,-2) rectangle (2,2);

  \node at (1.3,-1.3) {$W$};

  \node[centerdot] at (-1.1,0.4) {};
  \node[centerdot] at (-0.94,-0.45) {};
  \node[centerdot] at (-0.4,-1) {};
  \node[centerdot] at (0.3,-0.6) {};
  \node[centerdot] at (.73,0) {};

  \draw [-{Latex[length=3mm]}] (-0.76,-0.75) -- (-0.87,-0.61);

  % Outer Omega shape
  \draw[thick]
  (-1.5,2) .. controls (-1.8,0.5) and (-1.5,-1.6) .. (-0.5,-1.6)
   .. controls (0,-1.7) and (1,-1) .. (2,2);

  \node at (-0.5,-1.3) {$\Omega$};

% Omega(1) shape

  \clip
  (-1,2) .. controls (-1.3,0.4) and (-1,-1) .. (-0.4,-1)
   .. controls (0,-1.1) and (.6,-.6) .. (1.7,2);
  \fill[blue!20] (-2,-2) rectangle (2,2);

  \node[centerdot] at (-1.1,0.4) {};
  \node[centerdot] at (-0.94,-0.45) {};
  \node[centerdot] at (-0.4,-1) {};
  \node[centerdot] at (0.3,-0.6) {};
  \node[centerdot] at (.73,0) {};

  \draw [-{Latex[length=3mm]}] (-0.76,-0.75) -- (-0.87,-0.61);

  \draw[thick]
  (-1,2) .. controls (-1.3,0.4) and (-1,-1) .. (-0.4,-1)
   .. controls (0,-1.1) and (.6,-.6) .. (1.7,2);

  \node at (-0.25,-0.5) {$\Omega(1)$};

  % Inner Omega(2) region
  \clip
  (-0.5,2) .. controls (-0.8,0.4) and (-.7,-.4) .. (-0.2,-0.1)
   .. controls (0.3,0.2) and (.6,0.7) .. (1.2,2);
   \fill[red!20] (-2,-2) rectangle (2,2);
  \draw[thick]
  (-0.5,2) .. controls (-0.8,0.4) and (-.7,-.4) .. (-0.2,-0.1)
   .. controls (0.3,0.2) and (.6,0.7) .. (1.2,2);

   \node at (-0.1,.6) {$\Omega(2)$};

   % Arrow
   %\draw [-{Latex[length=3mm]}] (-0.5,-0.6) -- (-0.65,-0.45);

\end{scope}

% Arrow Psi to first D
\draw[->,thick] (0,-2.2) -- (0,-2.9) node[midway,left] {$\Psi$};

% First D with geodesics
\begin{scope}[shift={(0,-4.4)}]
  \fill[violet!20] (0,0) circle (1.2);
  \draw[thick] (0,0) circle (1.2);
  \node at (0,0) {$\mathbb D$};

  % Geodesics and points
  \draw[thick] (-1.2,0) arc (-96:70:.3);
  \draw[thick] (-1.2,0) arc (96:-70:.3);
  \draw[thick] ({1.2*cos(152)},{1.2*sin(152)}) arc (-120:50:.2);
  \draw[thick] ({1.2*cos(208)},{1.2*sin(208)}) arc (120:-50:.2);
  \draw[thick] ({1.2*cos(132)},{1.2*sin(132)}) arc (-140:30:.18);
  \draw[thick] ({1.2*cos(228)},{1.2*sin(228)}) arc (140:-30:.18);

  % Labels on arcs
  \node at ({1.44*cos(166)},{1.44*sin(166)}) {$K$};
  \node[centerdot] at ({0.91*cos(166)},{0.91*sin(166)}) {};
  \node at ({0.6*cos(166)},{0.6*sin(166)}) {$a_K$};

% Dots on |z|=1

  \node[centerdot] at (1.2,0) {};
  \node at (1.4,0) {$1$};
  \node[centerdot] at (-1.2,0) {};
  \node[centerdot] at ({1.2*cos(152)},{1.2*sin(152)}) {};
  \node[centerdot] at ({1.2*cos(208)},{1.2*sin(208)}) {};
  \node[centerdot] at ({1.2*cos(132)},{1.2*sin(132)}) {};
  \node[centerdot] at ({1.2*cos(228)},{1.2*sin(228)}) {};
  \node[centerdot] at ({1.2*cos(115)},{1.2*sin(115)}) {};
  \node[centerdot] at ({1.2*cos(245)},{1.2*sin(245)}) {};

  % Arrow on |z|=1
\draw [-{Latex[length=3mm]}] ({1.2*cos(195)},{1.2*sin(195)}) -- ({1.2*cos(205)},{1.2*sin(205)});

\end{scope}

% Arrow B to unit disk
\draw[->,thick] (0,-5.8) -- (0,-6.5) node[midway,left] {$B$};

% Second D: unit disk with |z|=1
\begin{scope}[shift={(0,-8)}]
  \fill[violet!20] (0,0) circle (1.2);
  \draw[thick] (0,0) circle (1.2);
  \node at (0,0) {$\mathbb D$};
  \node[centerdot] at (1.2,0) {};
  \node at (1.4,0) {$1$};

  % Arrow on |z|=1
\draw [-{Latex[length=3mm]}] (-1.2,0.1) -- (-1.2,-0.1);

\end{scope}

%---------------- RIGHT COLUMN: tau-plane -> exp -> annulus ----------------%

\draw[->,thick] (3.4,0) -- (4.1,0) node[midway,above] {$\tau$};

% Tau-plane with strip
\begin{scope}[shift={(8,0)}]

 % Right red region
  \fill[red!20] (-.8,-3) rectangle (2,2);
 % Blue strip
  \fill[blue!20] (-.8,-3) rectangle (0,2);

  % Axes
  \draw[->,thick] (-2.0,-.5) -- (2,-.5);
  \draw[->,thick] (-1.6,-3) -- (-1.6,2);
  %\node at (0,2.4) {$\tau$};

  % Vertical lines L1, L2
  \draw[thick] (-0.8,-3) -- (-0.8,2);
  \draw[thick] (0,-3) -- (0,2);

  % Dots on L1
  \node[centerdot] at (-0.8,-.5) {};
  \node[centerdot] at (-0.8,0.6) {};
  \node[centerdot] at (-0.8,-1.6) {};
  \node[centerdot] at (-0.8,1.7) {};
  \node[centerdot] at (-0.8,-2.7) {};
  \node[below left] at (-0.8,-0.5) {$1$};

% Arrow on L1
\draw [-{Latex[length=3mm]}] (-0.8,0.2) -- (-0.8,0.4);
  % Shaded strip between L1 and L2

  % Labels L1, L2, S
  \node[left] at (-.45,-3.3) {$L_1$};
  \node[left] at (0.35,-3.3) {$L_2$};
  %\fill (0,0) circle (1.5pt);
  \node[right] at (-.68,1.2) {$S$};
  \node[right] at (-1.5,1.2) {$\mathbb H_r$};
  \node[left] at (-0.75,-1.2) {$J$};

\end{scope}

% Arrow exp to annulus
\draw[->,thick] (8,-3.9) -- (8,-4.6) node[midway,right] {$\exp$};

% Annulus D(0,e)
\begin{scope}[shift={(8,-7)}]
  % Outer red region
  \fill[red!20] (-2,-2) rectangle (2,2);
  % Inner white disk
  %\fill[white] (0,0) circle (1);
  % Blue annulus (1 to e)
  \begin{scope}
    \clip (0,0) circle (2);
    \fill[blue!20] (0,0) circle (1.8);
    \fill[white] (0,0) circle (1);
  \end{scope}

  % Boundaries
  \draw[thick] (0,0) circle (1);
  \draw[thick] (0,0) circle (1.8);
  %\draw[thick] (0,0) circle (2);

  % Arrow on exp(L1)
\draw [-{Latex[length=3mm]}] (-1,0.1) -- (-1,-0.1);

  % Label
  \node at (0,0) {$D(0,e)$};
  \node[centerdot] at (1,0) {};
  \node at (1.2,0) {$e$};
\end{scope}

\end{tikzpicture}
\caption{\label{model} \small Key elements of the interpolation, including intervals~$J$ in the partition $\mathcal J$ of $L_1$ and the image arcs $K\in\mathcal K$ under $\Psi\circ \tau^{-1}$, with corresponding geodesics $\gamma_K$ in $\D$ and points $a_K$ on $\gamma_K$. The three colours are intended to show regions that correspond under the various maps.}
\end{figure}

Let $\mathcal J$ denote the regular partition of $L_1$ induced by the points $1+2\pi i \Z$ and let $\mathcal K$ denote the set of image arcs on $\partial\D$ of the intervals~$J$ in $\mathcal J$ under $\Psi\circ \tau^{-1}$. Since $\Psi$ and $\tau^{-1}$ are conformal mappings, adjacent arcs of the partition $\mathcal K$ have uniformly comparable lengths.

For any $K\in\mathcal K$, let $\gamma_K$ denote the hyperbolic geodesic in~$\D$ that joins the endpoints of~$K$ and let $a_K$ denote the point on $\gamma_K$ that is closest to~0. Then the Blaschke condition $\sum_{k\in \mathcal K}(1-|a_K|)<\infty$ holds, so the infinite product
\[
B(z)=\prod_{K\in\mathcal K}\frac{\ov{a_K}}{a_K}\frac{a_K-z}{1-\ov{a_K}z}
\]
is convergent, and defines a Blaschke product that maps~$\D$ onto~$\D$ and $\partial \D$ onto itself, the image curve winding round infinitely often.

At this point it might seem a good idea to attempt to interpolate across $\Omega(1)\setminus \ov{\Omega(2)}$ between
\begin{equation}\label{interpolate}
F(z) = \exp(\tau(z)),\;z\in\ov{\Omega(2)} \quad\text{and}\quad G(z) = e B\circ \Psi(z),\;z\in W,
\end{equation}
since~$F$ maps $\partial \Omega(2)$ to $\{z:|z|=e^2\}$ and~$G$ maps $\partial \Omega(1)$ to $\{z:|z|=e\}$, with their image curves winding round infinitely often in the same direction. However, there is no a priori reason why these image curves should wind round the circle at the same rate, which would allow such an interpolation to give a quasiregular mapping of~$\C$.

Bishop's idea here is to show (by ingenious use of the hyperbolic metric and harmonic measure estimates in $\D$) that there is a subset $\mathcal K_B\subset \mathcal K$ such that if~$B$ is defined instead as the infinite Blaschke product
\[
B(z)=\prod_{K\in\mathcal K_B}\frac{\ov{a_K}}{a_K}\frac{a_K-z}{1-\ov{a_K}z},
\]
then the sequence of points $e^{i\theta}$ satisfying $B(e^{i\theta})=1$, arranged in the order that the image under~$B$ winds round $\partial\D$, gives rise, under $\tau\circ \Psi^{-1}$, to a new partition $\mathcal L$ of $L_1$ such that
\begin{equation}\label{Mbound}
\text{each element of }\mathcal L \text{ meets at least two and at most } M \text{ elements of }\mathcal J,
\end{equation}
where~$M$ is a uniform bound independent of the geometry of~$\Omega$.

Note that each segment on $L_1$ of the new partition $\mathcal L$ corresponds under $B\circ \Psi\circ \tau^{-1}$ to a single circuit of $\partial \D$ starting and finishing at the point~1. Also note that the statement \eqref{Mbound} shows that there is a uniform comparability between the partitions~$\mathcal L$ and~$\mathcal J$ in terms of the relative frequency of occurrence of their intervals on $L_1$. In Bishop's words these partitions are `almost the same'. This makes it plausible that with this revised Blaschke product~$B$ we can interpolate across~$S$ between
\begin{equation}\label{interpolateS}
\exp(z),\;z\in L_2,\quad\text{and}\quad e B\circ \Psi\circ\tau^{-1}(z), \; z\in L_1,
\end{equation}
in such a way that the corresponding interpolation in \eqref{interpolate} across $\Omega(1)\setminus \ov{\Omega(2)}$  gives a quasiregular mapping on~$\C$.

The required interpolation on~$S$ is done by composing four maps (in the paper $\sigma,\psi_3,\psi_2,\psi_1$) that are mostly quasiconformal on~$S$. It would take us too far to give all the details but we shall describe here what is arguably the key map, $\psi_3$. This is the one that deals with the mismatch between the two partitions $\mathcal L$ and $\mathcal J$, that is, the fact that there may be up to~$M$ elements of the regular partition~$\mathcal J$ that meet any given element of $\mathcal L$. It is this map $\psi_3$ that accounts for the name `quasiconformal folding'.

First, each $K\in \mathcal L$ is approximated by a block consisting of an odd number of adjacent elements of $\mathcal J$ (at most~$M+1$ of them), whose union is called $I_K$ and whose lowest element is called $J_K$. The sets $I_K$ form a partition of the line $L_1$.

Then, for $K\in \mathcal L$, the rectangle $R_K=I_K\times [1,2]$ is introduced and also the cut rectangle $U_K=R_K \setminus X_K$, where $X_K$ is the closed segment connecting the upper left corner of $R_K$ with the centre of $R_K$; see Figure~3, borrowed from \cite{BishopJLMS}. The key quasiconformal mapping is given by this lemma.
\begin{lem}[Simple folding]\label{folding}
There is a quasiconformal mapping $\psi_3$ from $U_K$ onto $R_K$ so that the following hold:
\begin{itemize}
\item[(1)] $\psi_3$ is the identity on $\partial R_K\setminus L_1$;
\item[(2)] $\psi_3^{-1}$ extends continuously to $\partial R_K$;
\item[(3)] $\psi_3$ maps $I_K$ linearly onto $J_K$;
\item[(4)] $\psi_3$ maps `opposite sides' of $X_K$ linearly onto symmetric halves of $I_K\setminus J_K$ -- recall that $I_K\setminus J_K$ consists of an even number of elements of~$\mathcal J$;
\item[(5)] $\psi_3$ has bounded dilatation, the bound depending only on the constant~$M$ in \eqref{Mbound}.
\end{itemize}
\end{lem}
The name `folding' comes from the geometric behaviour of the inverse map $\psi_3^{-1}$, which stretches each of the three segments of $I_K$ by varying factors and then folds the result, mapping $J_K$ onto $I_K$ itself and the two halves of $I_K\setminus J_K$ onto $X_K$. Thus the map $\psi_3$ is more accurately an `unfolding'. Bishop calls this process `simple folding' because his other papers \cite{BishopActa,BishopPLMS} have more complicated models and require much more complicated folding to deal with them.
\begin{figure}
\centering
\includegraphics[width=0.4\textwidth]{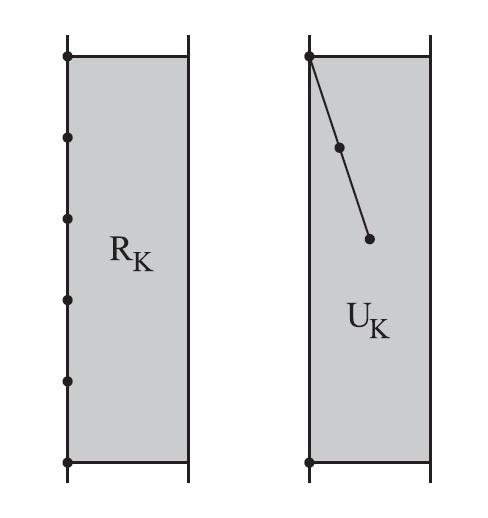}
\caption{\label{model} \small The rectangle $R_K$, with left edge consisting of five elements of the partition $\mathcal L$, and the cut rectangle $U_K$. The six dots on the left correspond under the folding map to the four on the right.}
\end{figure}

Finally, let~$U$ denote the domain~$S$ with all the segments $X_K$, $K\in\mathcal L$, removed. In view of properties~(1) and~(5), the collection of maps $\psi_3$ ranging over all such~$K$ give rise to a quasiconformal mapping of~$U$ onto~$S$, also called $\psi_3$. This maps each block $I_K$ of elements of~$\mathcal J$ onto $J_K$ and opposite sides of each cut $X_K$ onto $I_K\setminus J_K$. It is quasiconformal on~$U$ but not continuous (or even defined) on the cuts $X_K$. To fix this problem, a further map~$\sigma$  of~$S$ is applied after $\psi_3$, which ensures that the overall composed map $\sigma\circ\psi_3\circ\psi_2\circ\psi_1$ is continuous on the cuts $X_K$, quasiregular on~$S$, and interpolates in the way required by \eqref{interpolateS}. Hence the corresponding interpolated map in \eqref{interpolate} is quasiregular in~$\C$.

The proof of Theorem~\ref{Chris-jlms-main} is then completed using the Measurable Riemann Mapping Theorem, as described earlier.

\section{Further developments}

Before Bishop’s work, functions in the Eremenko--Lyubich class~$\mathcal B$ were constructed by first specifying a model $(\Omega, F)$ for the behaviour of the function near infinity, that is, a collection of unbounded simply connected domains, called tracts, each with a covering map to a punctured neighbourhood of infinity, as described in Section~1, and then using various bespoke approximation results. For example, in~\cite{RRRS} Rottenfusser et al used a method based on Cauchy integrals to construct functions in the class~$\mathcal B$ whose Julia sets are topologically far more complicated than was previously believed possible by using suitably `wiggly' tracts.

By contrast, there were few techniques before Bishop’s work for constructing functions in the Speiser class $\mathcal S$; for example, it was unclear whether the topologically complicated examples of Rottenfusser et al could arise in the class~$\mathcal S$. It was also unknown in general which models can be realised in the classes $\mathcal B$ and $\mathcal S$, and what, if any, differences exist between the dynamical behaviour of functions in the two classes.

Bishop’s paper \cite{BishopJLMS} greatly extends methods for constructing functions in the class~$\mathcal B$, in particular showing that every model consisting of unbounded simply connected domains can be realised. As applications, Rempe~\cite{arclike} used~\cite{BishopJLMS} to construct functions in the class $\mathcal B$ with Julia sets that contain every possible arc-like continuum, and Brown \cite{ABrown25} used it to construct functions in the class~$\mathcal B$ with `small' infinite order that are counterexamples to the so-called `strong Eremenko conjecture', that is, their escaping sets contain no curves to $\infty$. A somewhat different application can be found in \cite{BRSix} where the techniques in \cite{BishopJLMS} were adapted to construct an unbounded holomorphic function defined in the unit disc~$\D$ and in `class~$\mathcal B$ for the disc' such that~$f$ is bounded on a spiral that approaches $\partial \D$; in particular,~$f$ does not lie in the MacLane class~$\mathcal A$, consisting of holomorphic functions in $\D$ that have asymptotic values at points of a dense subset of $\partial \D$.

In two further papers, Bishop used more complicated versions of quasiconformal folding to construct an even wider range of functions in the classes~$\mathcal B$ and~$\mathcal S$, and so answer many open questions. In \cite{BishopActa}, his main theorem shows how to construct functions in the class~$\mathcal S$ with dynamical behaviour related to a given model $(\Omega,F)$ that again consists of tracts, unbounded simply connected domains, with associated covering maps. This time though the tracts of the model are defined to be the complementary components of an infinite planar tree satisfying certain mild geometric conditions. Moreover, the resulting entire function is not just in the class~$\mathcal S$ but it has no finite asymptotic values and exactly two critical values.

However, the conditions that the infinite tree must satisfy in \cite{BishopActa} are such that more tracts may need to be added to the original model for this construction to work, beyond those required to construct a function in the class~$\mathcal B$. The ideas of the construction in this paper link closely to the techniques developed in the paper \cite{BishopInv}, where Bishop solved Eremenko's question mentioned at the start of Section~1 of this review.

In \cite{BishopActa}, Bishop gives many striking applications of that paper's main theorem, perhaps the most notable being an answer to the longstanding open question as to whether functions in the class~$\mathcal B$ can have wandering domains. Recall from Section~1 that Eremenko and Lyubich had shown that functions in the class $\mathcal B$ cannot have escaping Fatou components, so any orbit of wandering domains of such functions must contain at least a subsequence accumulating in a bounded part of the plane. Bishop was able to develop his quasiconformal folding technique further in \cite{BishopActa}, replacing the infinite tree in the model by an infinite graph with both unbounded and bounded complementary components, each with an associated covering map onto one of the right half-plane, the left half-plane or the punctured unit disc. This allowed him to construct entire functions in the class~$\mathcal B$ with specified behaviour on certain bounded domains, and he thereby succeeded in constructing a function in the class~$\mathcal B$ with an oscillating wandering domain.

It remains a huge open question in transcendental dynamics whether there exists any entire function~$f$ with an orbit of wandering domains within which \textit{all} the limit values of the iterates $(f^n)$ lie in a bounded part of the plane.

In \cite{BishopPLMS}, Bishop proved further results related to the construction of class~$\mathcal S$ functions by quasiconformal folding, showing that for a given original model $(\Omega,F)$ there exists a function in class~$\mathcal S$ with no more than twice as many tracts as that model. Bishop also illuminates the differences between these two classes by showing that in the class $\mathcal S$ these additional tracts may sometimes be essential: for example, he shows that there is no function in the Speiser class whose only tract is a subdomain of a half-strip, whereas functions with this property are known to exist in the class $\mathcal B$.

Taken overall, Bishop’s three papers~\cite{BishopJLMS, BishopActa, BishopPLMS} on quasiconformal folding provide a `black box’ for constructing functions in the classes $\mathcal B$ and $\mathcal S$. For example, in~\cite{BishAlb}, Albrecht and Bishop used the main result on quasiconformal folding from \cite{BishopActa} to construct functions in the class~$\mathcal S$ for which the Hausdorff dimension of the Julia set takes values in the interval $(1,2]$ arbitrarily close to~1. Much earlier, Stallard had shown that for all functions in the class~$\mathcal B$ the Hausdorff dimension of the Julia set must lie in this interval, and moreover that the Hausdorff dimension of the Julia sets of such functions can take any value in this interval; see \cite{StallardII,StallardiV}.

Moreover, Bishop has adapted his original results on quasiconformal folding to solve (with co-authors) a variety of other major problems, not all in complex dynamics; see \cite{BishLaz19,BishLazUrb23,BishLaz24,BishRemp}, for example.

To conclude we mention some developments of Bishop's results by other authors.  In \cite{FGJ}, Fagella, Godillon and Jarque used the methods in a paper of Milhaljevi\'c-Brandt and Rempe-Gillen \cite{lasse-helena} to show that the class~$\mathcal B$ function with wandering domains built by Bishop in \cite{BishopActa} has only the wandering domains given by the construction. Then in \cite{FJL}, Fagella, Jarque and Lazebnik modified Bishop's construction to obtain a function in the class $\mathcal B$ with an orbit of wandering domains on all of which the function is univalent. Finally, we mention the paper by Mart\'i-Pete and Shishikura \cite{marshi} where a somewhat different method was used to construct functions in the class $\mathcal B$ with wandering domains and having finite order, in particular, any order that is a positive integer multiple of~$1/2$ (the lowest possible order in the class $\mathcal B$).

{\bf Acknowlegements}
The author is grateful to the referees, and also to Chris Bishop, Ian Short, Dave Sixsmith and Gwyneth Stallard, for very helpful comments on this review.

\bibliography{Wandering2}
\end{document}